\documentclass[a4paper,12pt]{article}

\usepackage{amsmath}
\usepackage{graphicx}
\usepackage{natbib}
\usepackage[top=1in, bottom=1in, left=1in, right=1in]{geometry}

\begin{document}

\title{A Mathematica Program for heat source function of 1D heat equation reconstruction by three types of data}
\author{Tomasz M. Lapinski \small{(84tomek@gmail.com)}, \\Sergey Leble \small{(leble@mif.pg.gda.pl)} \\Gdansk University of Technology, \\Faculty of Applied Physics and Mathematics,
\\{\small ul. Narutowicza 11/12, 80-952 Gdansk, Poland,}}
\maketitle

\begin{abstract}
	We solve an inverse problem for the one-dimensional heat diffusion equation. We reconstruct the heat source function for the three types of data: 1) single position point and different times, 2) constant time and uniformly distributed positions, 3) random position points and different times. First we demonstrate reconstruction using  simple inversion of discretized Kernel matrix. Then we apply Tikhonov regularization for two types  of the parameter of regularization estimation. The first one, which is in fact exemplary simulation, is based on minimization of the distance in C space of reconstructed function to the initial source function. Second rule is  known as Discrepancy principle. We generate the data from the chosen source function. In order to get some measure of accuracy of reconstruction we compare the result with the function from which data was generated.  We also deliver corresponding application in symbolic computation environment of Mathematica. The program has a lot of flexibility, it can perform reconstruction for much more general input then one considered in the paper.
\end{abstract}

\begin{keywords}
inverse problem, heat equation, source reconstruction, Tikhonov regularization
\end{keywords}

\maketitle

%%%%%		Section 1 - Introduction		%%%%%
%%%%%%%%%%%%%%%%%%%%%%%%%%%%

\section{Introduction}
	\indent There is a lot of attention to the topic of inverse problems for heat equation, e.g. in the recent papers concerning the subject \cite{comp}, \cite{flux} and \cite{ramm}. It is particularly interesting because it has a wide range of real world applications. Due to fact that its diffusion counterpart models so many phenomenas in physics and Geo science, solving corresponding inverse problem will always be in the demand. An example here would be problem of identifying the source of pollution in channels or rivers, even considered in 1D modeling. By measuring the density of pollution at some places and times in the river, one can find its source, see for example  \cite{poll}. We find that of high usefulness would be a computer program which could perform this kind of tasks. \\
	\indent In this article we  deliver a flexible symbolic computer program for Mathematica environment which performs the algorithm and all mentioned in the paper calculations.  The main task of the program  is to reconstruct the source function of the diffusion equation. The corresponding inverse problem, by Fourier series expansion is represented as Fredholm equation of the first kind with continuous kernel. It is well known that inverse operator of such integral operator is not continuous. Hence, the the solution of the equation does not depend continuously on the data in conventional Banach spaces, it is unstable. Therefore this is an ill-posed problem. The discretization of the equation results that the inverse operator is continuous, however, corresponding matrix representation has a very high condition number. Even small fluctuation in the data, for example the cutoff related to the precision of the computations makes reconstruction impossible.  We obtain so-called ill-conditioned problem. To obtain a well-posed problem we apply Tikhonov Regularization (see \cite{tih}, \cite{tih1}, \cite{tih2}). We consider two methods of choosing parameter of regularization. First one is by minimizing the difference between the original source function and the reconstructed one. For this method the knowledge of the source function is necessary, therefore we use this method only to check whether Tikhonov regularization gives satisfactory results. The second method, known as Discrepancy principle calculates the parameter using the error which the data is saddled with.  This method can be used to determine unknown source function, therefore it has a real world applications. The data used for computations is obtained by solving direct problem. \\
	\indent   The algorithm and the program  has a lot of places for adjustment, ex, number of data points and interval between them, length of the bar on which heat diffuses and many others. It is attached in the Appendix.

%%%%%		Section 2 - Direct Problem. Data generation		%%%%%
%%%%%%%%%%%%%%%%%%%%%%%%%%%%%%%%%%%%%

\section{Direct problem. Data generation}

	 We consider model given as differential equation at $x\in [0,l], t\in[0,\infty)  $
	\begin{equation}
		\label{introduction_differential_equation}
		\frac{\partial u(x,t)}{\partial t}=a^{2}\frac{\partial^{2}u(x,t)}{\partial x^{2}}+f(x),\\
	\end{equation}
	where $f(x)$ is a source function and $a$ a heat diffusion coefficient.\\
	We use Dirichlet boundary conditions
	\begin{equation}
		\label{introduction_boundary_conditions}
		u(0,t)=u(l,t)=0,
	\end{equation}
	and initial conditions
	\begin{equation}
		\label{introduction_intitial_conditions}
		u(x,0)=0,\notag
	\end{equation}
	where $l=1$ is length of bar.\\
	First we study principle features of the inverse problem via a simulation of data by means of this direct problem.
	The data for the inverse problem is obtained through evaluation of the   direct problem solution for a simple given pulse source function. We solve the direct problem and calculate the temperature for certain time and point on the  bar. The Fourier series method of  separation of variables is used to solve equation (\ref{introduction_differential_equation}) under conditions (\ref{introduction_boundary_conditions}) for $f(x)$ given by the function
	\begin{equation}
		\label{direct_problem_source_function}
		f(x)=\left\{
		\begin{array}{cc}
			0&[0,\frac{1}{3})\cup(\frac{2}{3},1] \\
			1&[\frac{1}{3},\frac{2}{3}]\\
		\end{array}
		\right.
	\end{equation}
	The solution of (\ref{introduction_differential_equation}) given (\ref{introduction_boundary_conditions}) and (\ref{direct_problem_source_function}) has following form
	\begin{equation}
		\label{direct_problem_solution}
		u(x,t)=\sum_{k=1}^{n} \sin(k \pi x)
    \bigg(\frac{1}{k \pi}\bigg)^{2}\Big(1- \exp\big[-(k \pi)^{2}t\big]\Big)\frac{1}{k \pi}\bigg(\cos\bigg(\frac{ k \pi }{3}\bigg)-\cos\bigg(\frac{2 k \pi }{3}\bigg)\bigg).
	\end{equation}
	with $a=1$ and for the simulation purposes the series has been cut to $n=100$ because increasing that number has no visible influence on the solution.\\
	We generate temperature data using solution (\ref{direct_problem_solution}). We do it on three ways
	\begin{enumerate}
	\item Data in one point $x_{0}$ but in different times $t_{0}$, starting with $t_{0}=0.004$ and separated by equal time interval $h_{t}$:
		\begin{equation*}
			u(x_{0},t_{i})=\sum_{k=1}^{n}\sin(k \pi x_{0})
    \bigg(\frac{1}{k \pi}\bigg)^{2}\Big(1-\exp\big[-(k \pi)^{2}t_{i}\big]\Big)\frac{1}{k \pi}\bigg(\cos\bigg(\frac{ k \pi }{3}\bigg)-\cos\bigg(\frac{2 k \pi }{3}\bigg)\bigg).
		\end{equation*}
   		where we consider $100$ points, i.e. $i=1,2,...,100$, separated by $h_{t}=0.004$, 
	\item Data at points $x_{i}$ on bar separated by equal length $h_{x}$ and constant time $t_{0}$:
 		\begin{equation*}
			u(x_{i},t_{0})=\sum_{k=1}^{n} \sin(k \pi x_{i})
    \bigg(\frac{1}{k \pi}\bigg)^{2}\Big(1-\exp\big[-(k \pi)^{2}t_{0}\big]\Big)\frac{1}{k \pi}\bigg(\cos\bigg(\frac{ k \pi }{3}\bigg)-\cos\bigg(\frac{2 k \pi }{3}\bigg)\bigg).
		\end{equation*}
		where we consider $100$ points separated by $h_{x}=\frac{1}{101}$, starting with $x_{1}=\frac{1}{101}$.
	\item Data in random points on bar $x_{r}$ and in different times $t_{i}$, starting with $t_{0}=0.004$ and separated by equal time interval $h_{t}$:
		\begin{equation*}
			u(x_{r},t_{i})=\sum_{k=1}^{n}\sin(k \pi x_{r})
    \bigg(\frac{1}{k \pi}\bigg)^{2}\Big(1-\exp\big[-(k \pi)^{2}t_{i}\big]\Big)\frac{1}{k \pi}\bigg(\cos\bigg(\frac{ k \pi }{3}\bigg)-\cos\bigg(\frac{2 k \pi }{3}\bigg)\bigg).
		\end{equation*}
		where $i=1,2,...,100$, $x_{r}$ is random value form $0$ to $1$.
	\end{enumerate}
	\begin{figure}[htbp]
		\begin{center}
			\includegraphics[scale=0.5]{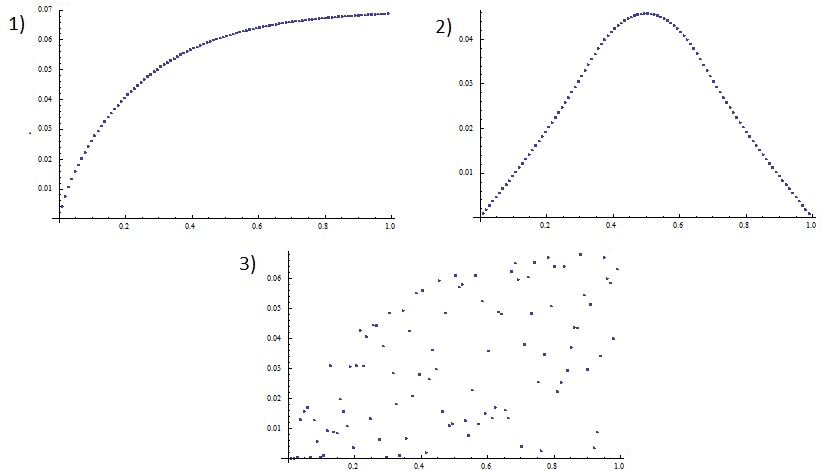}
			\caption{\em Simulation of three types of data 1) single $x_{0}$ and different times $t_{i}$, 2) constant time $t_{0}$ and uniformly distributed $x_{i}$, 3) random $x_{r}$ and different times $t_{i}$.}
		\end{center}
	\end{figure}
	Note, that values of data $u$ are saddled with the error
	\begin{equation}
		\label{direct_problem_delta}
		|u-u'|\leq\delta,
	\end{equation}
	where $u$ is data vector without error.\\
	We assume that the real data would be gathered with some device which has certain measurement error. In our case, the data is the temperature and is obtained by the ordinary thermometer has a measurement error $\delta=0.1$ but let us assume that we have gathered it with a special thermometer, very precise one and $\delta=0.0001$.
	The size of cutoff related to the precision of the computations made by program is much smaller, therefore does not affect $\delta$. \\
	We include the measurement error in the data by adding to it a random component of order of $\delta$
	\begin{equation}
		\label{direct_problem_measurement_error}
		u'=u+\delta r
	\end{equation}
	where $r$ is some constant random vector such that (\ref{direct_problem_delta}) is valid.

%%%%%		Section 3 - Inverse problems. Methods		%%%%%
%%%%%%%%%%%%%%%%%%%%%%%%%%%%%%%%%%

\section{Inverse problem. Methods}

%%%%%		Subsection 1 - Naive Method		%%%%%
	\subsection{Naive Method}
		Solution (\ref{direct_problem_solution}) can by represented in the form of Fredholm equation of the first kind as
		\begin{equation}
			\label{inverse_problem_fredholm_equation}
			u((x,t)_{i})=\int_{0}^{1}K((x,t)_{i},s)f(s)ds.
		\end{equation}
		where $u$ is temperature data $(x,t)_{i}$ point of data and for the three types of data, respectively, Kernel has the form 
		\begin{align}
			\label{inverse_problem_kernel}
			&K((x,t)_{i},s)=\sum_{k=1}^{n} \sin( k \pi x_{i})
    \bigg(\frac{1}{ k \pi}\bigg)^{2}\Big(1- \exp\big[-( k \pi)^{2}t_{i}\big]\Big)\sin(k \pi s),\\
			&\ 1. \ t_{i}=t_{0}, \notag\\
			&\ 2. \ x_{i}=x_{0}, \notag\\
			&\ 3. \ x_{i}=x_{r},\notag
		\end{align}
		After the discretization, equation (\ref{inverse_problem_fredholm_equation}), including Kernel (\ref{inverse_problem_kernel}) is given by
		\begin{align}
			\label{inverse_problem_discretized_fredholm_equation}
			&u((x,t)_{i})=\sum_{k=1}^{n} \sin( k \pi x_{i})
    \bigg(\frac{1}{ k \pi}\bigg)^{2}\big(1- \exp\big[-( k \pi)^{2}t_{i}\big]\big)\sum_{j=1}^{J}\sin(k \pi s_{j})f(s_{j})h_{s},\notag\\
			&\ 1. \ t=t_{0},\\
			&\ 2. \ x=x_{0},\notag\\
			&\ 3. \ x=x_{r}.\notag
		\end{align}
		where we consider $J=100$, which is a number of points of heat source function in the equation.\\ 
		Equation (\ref{inverse_problem_discretized_fredholm_equation}) can be represented as matrix equation
		\begin{equation}
			u=Kfh_{s}.
		\end{equation}
		where $h_{s}=\frac{1}{101}$ are interval between points. \\
		To obtain source matrix $f$, we simply invert the matrix $K$.
		\begin{equation*}
 			f=\frac{1}{h_{s}}K^{-1}u'.
		\end{equation*}
		where we use the data saddled with the measurement error given by (\ref{direct_problem_measurement_error}).\\
		After performing computations using programming environment Mathematica, it appears that, with this accuracy (100 measurement), for all three types of the data method gives results which significantly differs from original source function (see Figure 2). For the data type 1) and 3) the program could not calculate the inverse properly of the matrix becaues of its high condition number. The reason this results is the ill-conditioning of the problem. Since the considered matrix has a very high condition number, the error related to the data $u'$ which for our program is of range $\delta=0.0001$ significantly affects the solution $f$. 
		\begin{figure}[h]
			\begin{center}
				\includegraphics[scale=0.5]{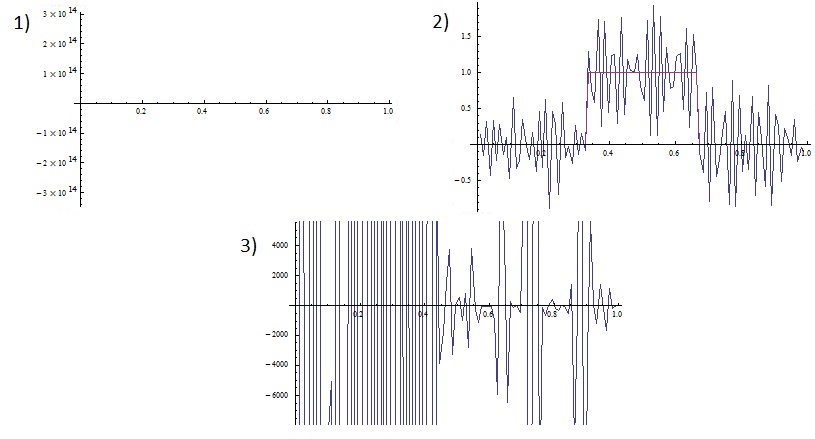}\\
				\caption{\em Results of reconstruction for the naive method for three types of data, 1) single $x_{0}$ and different times $t_{i}$, 2) constant time $t_{0}$ and uniformly distributed $x_{i}$, 3) random $x_{r}$ and different times $t_{i}$.}
			\end{center}
		\end{figure}

%%%%%		Subsection 2 - Tikhonov regularization		%%%%%
	\subsection{Tikhonov Regularization}
		We obtain a well-posed problem by applying Tikhonov Regularization (see \cite{tih}, \cite{tih1}, \cite{tih2}, \cite{regularization_of_inverse_problems}) to original problem. This method is based on Tikhonov Functional which for our case has following form
 		\begin{equation}
			\label{inverse_problems_tikhonov_functional}
			M^{\alpha}(f)=\int_{a}^{b}\bigg[\int_{0}^{1}K(x_{i},t_{i},s)f(s)ds-u(x_{i},t_{i})\bigg]^{2}dx+\alpha\int_{0}^{1}f^{2}(s)ds,\\
		\end{equation}
		It consist of two parts. First, is ordinary least square part and the second is norm minimization part.
		For different types of data we have different variables
		\begin{align*}
			&1. \ t_{i}=t_{0},\ a=0.004,\ b=N*h_{t}+a\\
			&2. \ x_{i}=x_{0},\ a=\frac{1}{101},\ b=\frac{100}{101},\\
			&3. \ x_{r}=x_{i},\ a=0.004,\ b=N*h_{t}+a.
		\end{align*}
		Discretized functional (\ref{inverse_problems_tikhonov_functional}) has a following form
		\begin{equation*}
			M^{\alpha}(f)=\sum_{i=1}^{N}h_{i}\bigg[\sum_{j=1}^{J}h_{s} K_{ij}f_{j}-u_{i}\bigg]^{2}+\alpha\sum_{j=1}^{J} h_{s} f_{j}^{2}.
		\end{equation*}
		where $\alpha$ is a parameter of regularization, $J=100$ is a number of heat source function and $N=100$ is number of temperature data points.\\
		We minimize discretized functional by differentiation w.r.t. $f_{j}$ and compare to $0$. We obtain 
		\begin{equation*}
			\sum_{i=1}^{N}\sum_{j=1}^{J}h_{i}h_{s} K_{ij} K_{ik}  f_{j}+\alpha h_{s} f_{k}-\sum_{i=1}^{N} h_{i} K_{ik}u_{i}=0.
		\end{equation*}
		\indent After changing representation from discrete to matrix, the solution has form
		\begin{equation}
			\label{inverse_problems_tikhonov_formula}
			f_{\alpha}=(K^{T}K h_{i}h_{s} +\alpha h_{s} I)^{-1}K^{T}  h_{i} u'.
		\end{equation}
		wher $I$ is identity matrix and we use the data saddled with the measurement error given by (\ref{direct_problem_measurement_error}) and for each data type we have
		\begin{align*}
			&1. \ h_{i}=h_{x},\\
			&2. \ h_{i}=h_{t},\\
			&3. \ h_{i}=h_{t},\
		\end{align*}
		To compute unknown source function $f$, knowledge of a value of parameter $\alpha$ is needed.\\
		We consider following parameter choice rules.

%		Subsubsection 1 - Minimization of the distance to initial data		%
		\subsubsection{Minimization distance to initial data}
			First rule is finding $\alpha$ by minimizing the square error between source function used in simulated measurement and source function obtained from calculation of the inverse problem, i.e.
			\begin{equation}
				\label{inverse_problems_tikhonov_alpha_minimazing_distance}
				\alpha=\min \sum_{i=1}^{100}\bigg[f_{i}-\big[(K^{T}K h_{i}h_{s} +\alpha h_{s} I)^{-1}K^{T}h_{i}u'\big]_{i}\bigg]^{2}
			\end{equation}
			where $f$ here is a discretized initial source function (\ref{direct_problem_source_function}). The reconstructed source function for obtained this way parameter $\alpha$ is shown on Figure 3.
			\begin{figure}[h]
				\begin{center}
					\includegraphics[scale=0.5]{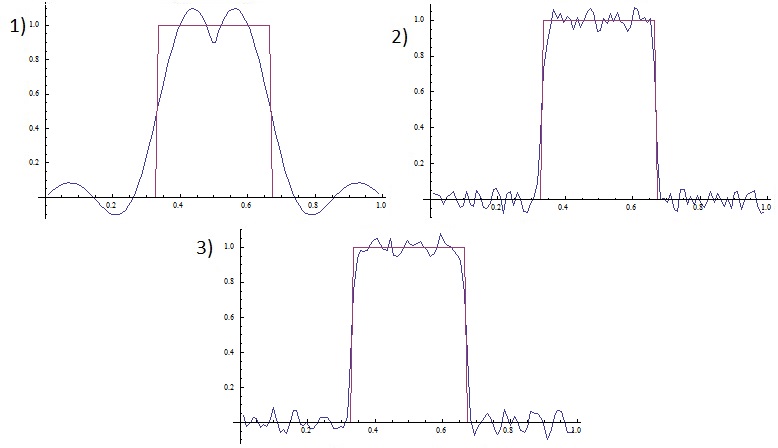}
					\caption{\em Results of reconstruction of Tikhonov method for three types of data with choosing $\alpha$ by minimizing distance of outcome to initial data, 1) single $x_{0}$ and different times $t_{i}$, 2) constant time $t_{0}$ and uniformly distributed $x_{i}$, 3) random $x_{r}$ and different times $t_{i}$.}
				\end{center}
			\end{figure}
			The algorithm which find the minimum start calculations with certain large $\alpha$ and decreases it when sum in (\ref{inverse_problems_tikhonov_alpha_minimazing_distance}) decreases by multiplying with some constant. Program repeats the step of decreasing $\alpha$ until the difference does not decrease further.

%		Subsubsection- Dicrepancy principle		%
		\subsubsection{Discrepancy principle}
			The regularization parameter defined by discrepancy principle (for details see \cite{regularization_of_inverse_problems}) is
			\begin{equation}
				\label{inverse_problems_tikhonov_alpha_discrepancy_principle}
				\alpha=\sup\{\alpha>0\ |\ \|Kf_{\alpha}h_{x}-u'\|<\delta\},
			\end{equation}
			where $f_{\alpha}$ is a reconstructed source function using equation (\ref{inverse_problems_tikhonov_formula}), $u'$ is used data, $h_{x}$ is an interval between points of discretization and $\delta$ is an measurement error of the data, but in our particular case it is precision of the computations. The outcomes for three types of data using this parameter choice rules are  given by Figure 4.
			\begin{figure}[h]
				\begin{center}
					\includegraphics[scale=0.5]{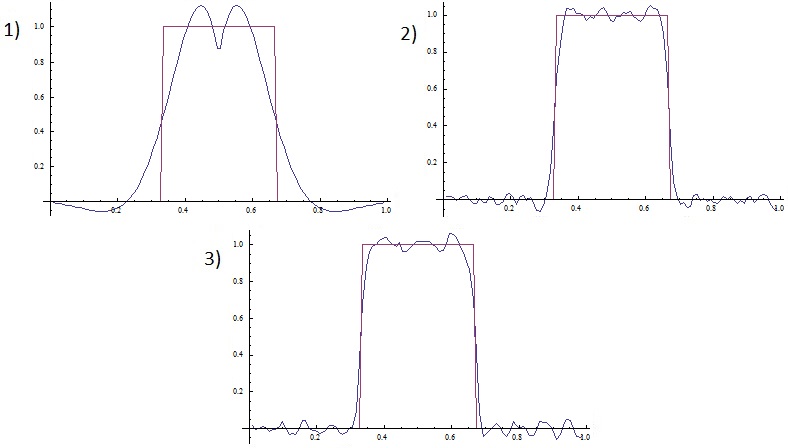}
					\caption{\em Results of reconstruction of Tikhonov method for three types of data with $\alpha$ chosen by Discrepancy Principle, 1) single $x_{0}$ and different times $t_{i}$, 2) constant time $t_{0}$ and uniformly distributed $x_{i}$, 3) random $x_{r}$ and different times $t_{i}$.}
				\end{center}
			\end{figure}
			The algorithm start with some very small $\alpha$ and increases it by multiplying with some constant when the condition in supremum of (\ref{inverse_problems_tikhonov_alpha_discrepancy_principle}) is still valid. It repeats the step until the condition is no longer valid.

%%%%%		Section 4 - Conclusions		%%%%%
%%%%%%%%%%%%%%%%%%%%%%%%%%%%

\section{Conclusion}
	In this paper we first show that without regularization the solution is far from desired, as expected, because of instability. Using regularization, the reconstructed functions are very close to original ones, for both type of parameter choice rules. From obtained results we can conclude that application of Tikhonov Regularization is very effective in considered problems.

%%%%%		Section 5 - Appendix - Program		%%%%%
%%%%%%%%%%%%%%%%%%%%%%%%%%%%%%%

\appendix
\section{Appendix - Program}

	\includegraphics{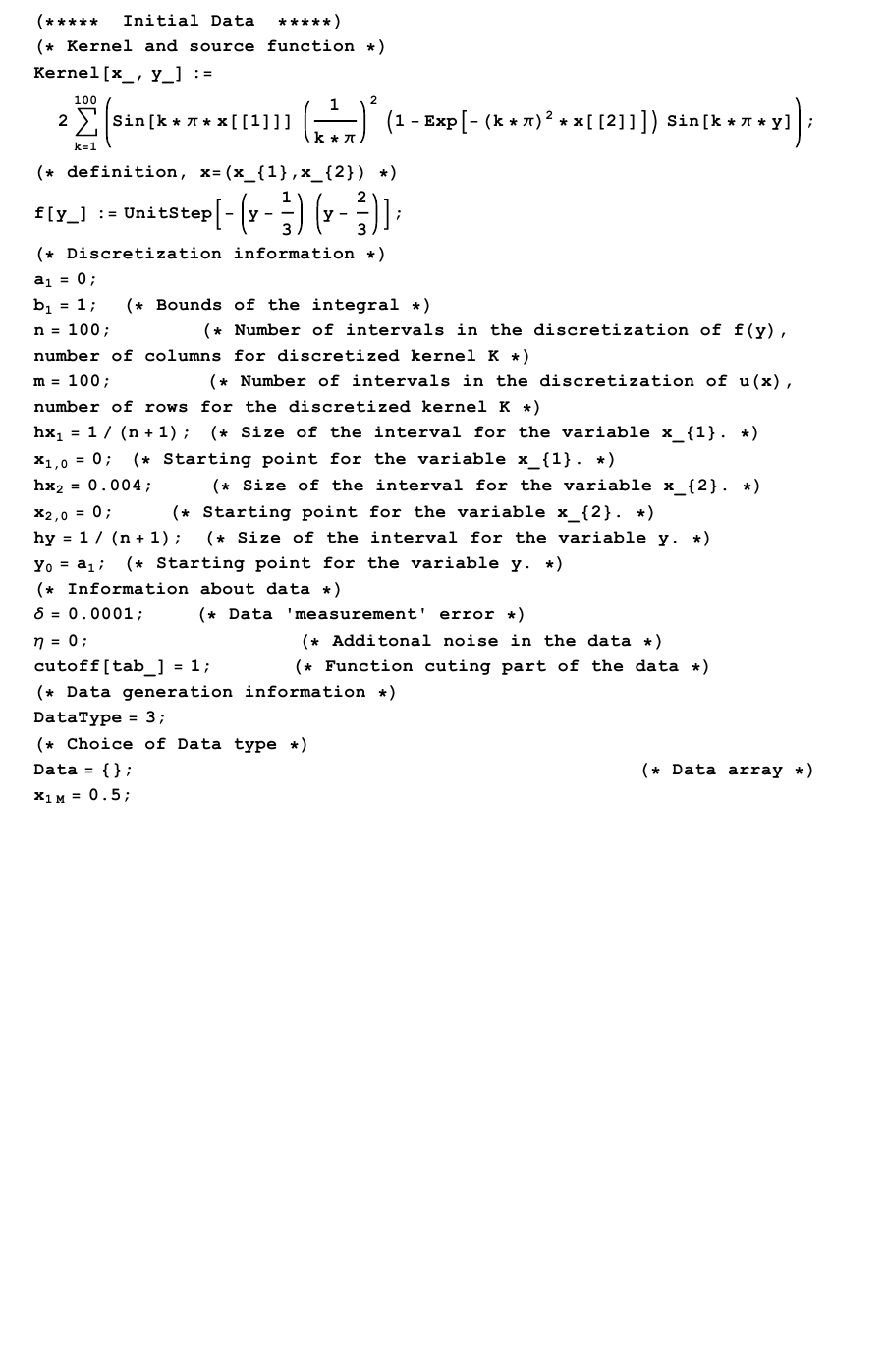}\\
\\
	\includegraphics{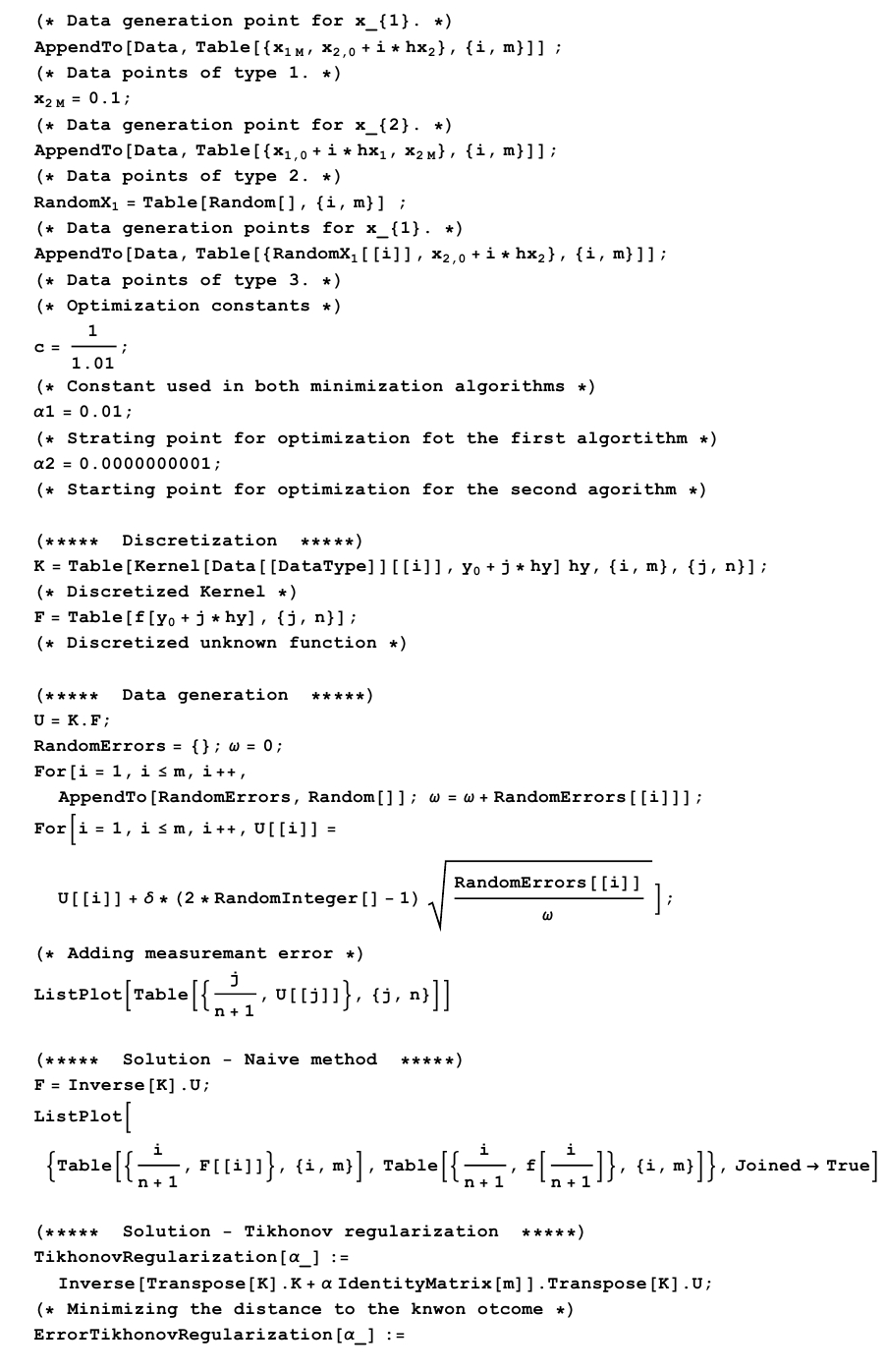}\\
	\includegraphics{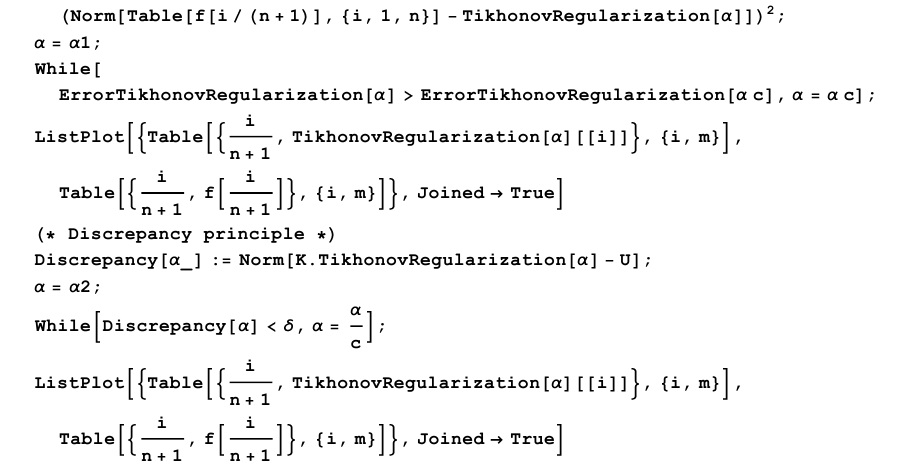}

%%%%%		Bibliography		%%%%%
%%%%%%%%%%%%%%%%%%%%%%

\bibliographystyle{plainnat}
\bibliography{biblio}

\end{document}